\documentclass{amsproc}
\usepackage{amsfonts}
\usepackage{makecell}
\usepackage{lscape}

\theoremstyle{plain}

\newtheorem{corollary}{Corollary}

\newtheorem{definition}{Definition}

\newtheorem{proposition}{Proposition}
\newtheorem{remark}{Remark}

\newtheorem{theorem}{Theorem}
\newcommand{\mf}[1]{\mathfrak{#1}}

\numberwithin{equation}{section}

\begin{document}
\title[Absolute continuity]{The Absolute continuity of convolution products
of orbital measures in exceptional symmetric spaces}
\author{K. E. Hare}
\address{Dept. of Pure Mathematics\\
University of Waterloo\\
Waterloo, Ont.,~Canada\\
N2L 3G1}
\email{kehare@uwaterloo.ca}
\thanks{This research is supported in part by NSERC \#44597.}
\author{Jimmy He}
\address{Dept. of Pure Mathematics\\
University of Waterloo\\
Waterloo, Ont.,~Canada\\
N2L 3G1}
\email{jimmy.he@uwaterloo.ca}
\thanks{This paper is in final form and no version of it will be submitted
for publication elsewhere.}

\begin{abstract}
Let $G$ be a non-compact group, let $K$ be the compact subgroup fixed by a
Cartan involution and assume $G/K$ is an exceptional, symmetric space, one
of Cartan type $E,F$ or $G$. We find the minimal integer, $L(G),$ such that
any convolution product of $L(G)$ continuous, $K$-bi-invariant measures on $%
G $ is absolutely continuous with respect to Haar measure. Further, any
product of $L(G)$ double cosets has non-empty interior. The number $L(G)$ is
either $2$ or $3$, depending on the Cartan type, and in most cases is
strictly less than the rank of $G$.
\end{abstract}

\maketitle

\section{\protect\bigskip Introduction}

In \cite{Ra} Ragozin showed that if $G$ is a connected simple Lie group and $%
K$ a compact, connected subgroup fixed by a Cartan involution of $G$, then
the convolution of any $\dim G/K$, continuous, $K$-bi-invariant measures on $%
G$ is absolutely continuous with respect to the Haar measure on $G$. This
was improved by Graczyk and Sawyer in \cite{GSJFA}, who showed that when $G$
is non-compact and $n=$ rank $G/K$, then the convolution product of any $n+1$
such measures is absolutely continuous. They conjectured that $n+1$ was
sharp with this property. This conjecture was shown to be false, in general,
for the symmetric spaces of classical Cartan type in \cite{GHSym} (where it
was shown that rank $G/K$ suffices except when the restricted root system is
type $A_{n}$). Here we prove the conjecture is also false for the symmetric
spaces of exceptional Cartan type (Types $EI-IX,$ $FI,$ $FII$ and $G),$
except when the restricted root system is type $BC_{1}$ or $A_{2}$. In fact,
the sharp answer for the exceptional symmetric spaces (as detailed in Cor. %
\ref{maincor}) is always either two or three, depending on the type, even
though the rank can be as much as 8.

A special example of $K$-bi-invariant continuous measures on $G$ are the 
\textit{orbital measures }%
\begin{equation*}
\nu _{z}=m_{K}\ast \delta _{z}\ast m_{K},
\end{equation*}%
where $m_{K}$ is the Haar measure on $K$. These singular probability
measures are supported on the double cosets $KzK$ in $G$. One consequence of
our result is that any product of three (and often two) double cosets has
non-empty interior in $G$. Questions about the convolution of orbital
measures are also of interest because of the connection with the product of
spherical functions, c.f. \cite{GS2003}, \cite{GS2004}.

To establish our result, we actually study an equivalent absolute continuity
problem for measures on $\mathfrak{g}$, the Lie algebra of $G$. The Cartan
involution gives rise to a decomposition of the Lie algebra as $\mathfrak{%
g=k\oplus p}$ where $\mathfrak{k}$ is the Lie algebra of $K$. By the $K$%
-invariant measures on $\mathfrak{p,}$ we mean those satisfying $\mu (E)=\mu
(Ad(k)E)$ for all $k\in K$ and Borel sets $E$. An example of a $K$%
-invariant, singular probability measure on $\mathfrak{p}$ is the orbital
measure $\mu _{Z}$ defined for $Z\in \mathfrak{p}$ by 
\begin{equation*}
\int_{\mathfrak{p}}fd\mu _{Z}=\int_{K}f(Ad(k)Z)dm_{K}(k)
\end{equation*}%
for any continuous, compactly supported function $f$ on $\mathfrak{p}$. This
measure is supported on the $Ad(K)$-orbit of $Z$, denoted

\begin{equation*}
O_{Z}=\{Ad(k)Z:k\in K\}\text{.}
\end{equation*}

It was shown in \cite{AG} that $\mu _{Z_{1}}\ast \cdot \cdot \cdot \ast \mu
_{Z_{L}}$ is absolutely continuous with respect to Lebesgue measure on $%
\mathfrak{p}$ if and only if $\nu _{z_{1}}\ast \cdot \cdot \cdot \ast \nu
_{z_{L}}$ is absolutely continuous on $G$ when $z_{j}=\exp Z_{j},$ and the
former is the problem we will actually study.

A compact, connected simple Lie group $G$ can be viewed as the symmetric
space $(G\times G)/G$. The absolute continuity problem has been studied for
orbital measures in this setting, as well, with the sharp answers found for
the exceptional Lie groups in \cite{HJSY}. The arguments used there (as well
as for the classical Lie groups/algebras in \cite{GHAdv}, \cite{GHMathZ})
were based upon harmonic analysis methods not generally available in
symmetric spaces. Instead, here we rely heavily upon a combinatorial
condition for absolute continuity that was discovered by Wright in \cite{Wr}
for compact Lie algebras, and was extended to symmetric spaces in \cite%
{GHSym}. Another key idea in our work is an embedding argument based on the
Freudenthal magic square construction of particular exceptional symmetric
spaces. This embedding idea was inspired, in part, by an embedding argument
used in \cite{GSColloq} for particular classical symmetric spaces.

In \cite{RS1} Ricci and Stein proved that if a convolution product of
orbital measures is absolutely continuous, then its density function is
actually in $L^{1+\varepsilon }$ for some $\varepsilon >0$. It would be
interesting to know the size of $\varepsilon $, however our arguments do not
yield any information about this. For other approaches to the study of sums
of adjoint orbits, we refer the reader to \cite{DRW} and \cite{KT}, for
example.

\section{Notation and Statement of main results}

\subsection{Notation and Terminology}

Let $G$ be a non-compact, simple, connected Lie group with Lie algebra $%
\mathfrak{g}$ and let $\Theta $ be an involution of $G$. We assume that $%
K=\{g\in G:\Theta (g)=g\}$ is compact and connected. The quotient space, $%
G/K $, is called a symmetric space.

The map $\Theta $ induces an involution of $\mathfrak{g}$, denoted $\theta $%
. We put 
\begin{equation*}
\mathfrak{k}=\left\{ X\in \mathfrak{g\mid \theta }\left( X\right) =X\right\} 
\text{ and }\mathfrak{p}=\left\{ X\in \mathfrak{g\mid \theta }\left(
X\right) =-X\right\} ,
\end{equation*}%
the $\pm 1$ eigenspaces of $\theta $, respectively. The decomposition $%
\mathfrak{g}=\mathfrak{k\oplus p}$ is called the Cartan decomposition of the
Lie algebra $\mathfrak{g}.$ We fix a maximal abelian (as a subalgebra of $%
\mathfrak{g}$) subspace $\mathfrak{a}$ of $\mathfrak{p}$. It is known that $%
Ad(k):\mathfrak{p}\rightarrow \mathfrak{p}$ whenever $k\in K$ where $%
Ad(\cdot )$ denotes the adjoint action of $G$ on $\mathfrak{g}$. If we put $%
A=\exp \mathfrak{a}$, then $G=KAK$, hence every double coset contains an
element of $A$. Similarly, every $Ad(K)$-orbit contains an element of $%
\mathfrak{a,}$ so in studying the orbital measures $\nu _{x}$ or $\mu _{X}$
there is no loss in assuming $x\in A$ and $X\in \mathfrak{a}$.

By the absolute root system, we mean the root system of $\mathfrak{g}$. By
the restricted roots of the symmetric space we mean the restrictions to $%
\mathfrak{a}$ of the roots that are non-zero on $\mathfrak{a}$. This set,
called the restricted root system, will be denoted $\Phi $. The vector spaces%
\begin{equation*}
\mathfrak{g}_{\alpha }=\left\{ X\in \mathfrak{g:}\left[ H,X\right] =\alpha
\left( H\right) X\text{ for all }H\in \mathfrak{a}\right\} ,
\end{equation*}%
for $\alpha \in \Phi ,$ are known as the restricted root spaces. We remark
that these need not be one dimensional and that $\theta (\mathfrak{g}%
_{\alpha })=\mathfrak{g}_{-\alpha }$. The space $\mathfrak{p}$ can also be
described as 
\begin{equation*}
\mathfrak{p}={sp}\left\{ X-\theta X:X\in \mathfrak{g}_{\alpha },\ \alpha \in
\Phi ^{+}\right\} \oplus \mathfrak{a}
\end{equation*}%
where by $sp$ we mean the real span. We put 
\begin{equation*}
X^{+}=X+\theta X\text{ and }X^{-}=X-\theta X.
\end{equation*}%
Of course, $X^{+}\in \mathfrak{k}$.

Throughout the paper, we assume $G/K$ is an irreducible, Riemannian,
globally symmetric space of Type III, whose absolute root system is
exceptional. (We call these `exceptional symmetric spaces'). They have
Cartan classifications $EI-IX$, $FI,$ $FII$ and $G$. In the appendix, we
summarize basic information about these symmetric spaces, including the
absolute and restricted root systems, and the rank and dimension of $G/K$.
This information is taken from \cite{Bu} and \cite{He}.

Given $Z\in \mathfrak{a}$, we let 
\begin{equation*}
\Phi _{Z}=\{\alpha \in \Phi :\alpha (Z)=0\}
\end{equation*}%
be the set of \textit{annihilating (restricted) roots} of $Z$ and let 
\begin{equation*}
\mathcal{N}_{Z}=sp\{X_{\alpha }^{-}:X_{\alpha }\in \mathfrak{g}_{\alpha }%
\text{, }\alpha \notin \Phi _{Z}\}\subseteq \mathfrak{p}\text{.}
\end{equation*}%
The set $\Phi _{Z}$ is itself a root system that is $\mathbb{R}$\textit{%
-closed} (meaning $sp\Phi _{Z}\cap \Phi =\Phi _{Z})$ and is a proper root
subsystem provided $Z\neq 0$. By the \textit{type }of $Z$ we mean the Lie
type of $\Phi _{Z}$.

The annihilating root subsystems$,$ $\Phi _{Z}$, and the associated spaces, $%
\mathcal{N}_{Z},$ are of fundamental importance in studying orbits and
orbital measures. Indeed, the tangent space to the $Ad(K)$-orbit of $Z$ is $%
T_{Z}(O_{Z})=\mathcal{N}_{Z}$. In particular, $\dim O_{Z}=\dim \mathcal{N}%
_{Z},$ so $O_{Z}$ is a proper submanifold of $\mathfrak{p}$ and hence of
measure zero. It has positive dimension if $Z\neq 0$, in which case we call
the orbit $O_{Z}$ and the orbital measure $\mu _{Z}$ non-trivial. More
generally, if $X\in O_{Z}$, say $X=Ad(k)Z$, then $%
T_{X}(O_{Z})=Ad(k)T_{Z}(O_{Z})=Ad(k)\mathcal{N}_{Z}$.

Similar statements hold for the orbital measures on $G$ and double cosets.
In particular, we call an orbital measure $\nu _{z}$ non-trivial if $z=\exp
Z $ where $\Phi _{Z}\neq \Phi $. Equivalently, $z$ is not in the normalizer
of $K$ in $G$.

\subsection{Main Result}

Here is the statement of our main result, whose proof will be the content of
the remainder of the paper. Note that when we say an orbital measure is
absolutely continuous (or singular) we mean either with respect to Haar
measure on $G$ or Lebesgue measure on $\mathfrak{p}$, depending on whether
the orbital measure is on $G$ or $\mathfrak{p}.$

\begin{theorem}
\label{main} Let $G/K$ be an exceptional symmetric space. For each non-zero $%
X\in \mathfrak{a}$ there is an integer $L=L_{X}$ such that $\mu _{X}^{L}$
and $\nu _{\exp X}^{L}$ is absolutely continuous if $L\geq L_{X}$ and
singular otherwise. Moreover, $L_{X}=2$ other than in the following cases
when $L_{X}=3$:

\medskip

\begin{tabular}{|p{3.2cm}|c|c|c|c|c|}
\hline
Cartan type of $G/K$ & $EI$ & $EIV$ & $EV$ & $EVII$ & $EVIII$ \\ \hline
Type of $X$ & $D_5,A_5$ & $A_1$ & $E_6,D_6$ & $C_2$ & $E_7$ \\ \hline
\end{tabular}

\medskip

Let $L(G)=3$ if $G/K$ is one of type $EI$, $EIV$, $EV$, $EVII$ and $EVIII$
and let $L(G)=2$ otherwise. Then any convolution product of $L(G)$
non-trivial orbital measures on $G$ or $\mathfrak{p}$ is absolutely
continuous.
\end{theorem}

We continue to use the same definition of $L(G)$ in the corollary.

\begin{corollary}
\label{maincor}(1) Any product of $L(G),$ non-trivial double cosets and any
sum of $L(G),$ non-trivial $Ad(K)$-orbits has non-empty interior in $G$ and $%
\mathfrak{p}$ respectively.

(2) Any convolution product of $L(G),$ $K$-bi-invariant, continuous measures
on $G$ or $K$-invariant, continuous measures on $\mathfrak{p}$ is absolutely
continuous and this is sharp.
\end{corollary}

\begin{proof}
The first statement follows from the fact that $\mu _{Z_{1}}\ast \cdot \cdot
\cdot \ast \mu _{Z_{t}}$ is absolutely continuous with respect to Lebesgue
measure on $\mathfrak{p}$ if and only if the sum $O_{Z_{1}}+\cdot \cdot
\cdot +O_{Z_{t}}$ has non-empty interior if and only if $K\exp Z_{1}K\cdot
\cdot \cdot K\exp Z_{t}K$ has non-empty interior, see \cite{GHSym}.

Furthermore, in \cite{Ra} it was shown that if every product of $L$
non-trivial double cosets has non-empty interior, then any convolution
product of $L$ $K$-bi-invariant, continuous measures on $G$ is absolutely
continuous. By \cite{AG}, this is equivalent to any product of $L$ $K$%
-invariant continuous measures on $\mathfrak{p}$ being absolutely continuous.

The sharpness of $L(G)$ is immediate from the theorem.
\end{proof}

\begin{remark}
We note that the ranks of the symmetric spaces for which $L(G)=2$ vary from
1-4, while those with $L(G)=3$ range from 2-8.
\end{remark}

In section 3, we will prove the absolute continuity for $L\geq L_{X}$. In
section four we prove singularity for $L<L_{X}$. It will suffice to work
with orbital measures on $\mathfrak{p}$ because of the equivalence result of
Anchouche and Gupta \cite{AG}.

\section{Proof of Absolute continuity}

\subsection{Criteria for Absolute continuity}

There are two main ideas we use to establish absolute continuity. The first
is a combinatorial argument based on the annihilating root systems.

Note that by the \textit{rank of a root system} we mean the dimension of the
Euclidean space it spans. An annihilating root system $\Phi _{X},$ for $%
X\neq 0,$ always has proper rank. By the\textit{\ dimension of a root system}
$\Phi _{0}$ we mean 
\begin{equation*}
\dim \Phi _{0}=\dim sp\{E_{\alpha }^{-}:E_{\alpha }\in \mathfrak{g}_{\alpha }%
\text{, }\alpha \in \Phi _{0}\},
\end{equation*}%
that is, the cardinality of $\Phi _{0}^{+}$ counted by dimension of the
corresponding restricted root spaces, equivalently, the sum of the
multiplicities of the restricted positive roots in $\Phi _{0}$.

\begin{theorem}
\label{WrCriteria} \cite{GHSym}, \cite{Wr} Assume $G/K$ is a symmetric space
with restricted root system $\Phi $ and Weyl group $W$. Suppose $%
X_{1},...,X_{m}$ $\in \mathfrak{a}$ and assume%
\begin{equation}
(m-1)\left( \dim \Phi -\dim \Psi \right) -1\geq \sum_{i=1}^{m}\left( \dim
\Phi _{X_{i}}-\min_{\sigma \in W}\dim (\Phi _{X_{i}}\cap \sigma (\Psi
))\right)  \label{WrC}
\end{equation}%
for all $\mathbb{R}$-closed, root subsystems $\Psi \subseteq \Phi $ of
co-rank $1$. Then $\mu _{X_{1}}\ast \cdot \cdot \cdot \ast \mu _{X_{m}}$ is
absolutely continuous.
\end{theorem}

An easy calculation shows the following.

\begin{corollary}
Suppose $X_{1},...,X_{m}$ $\in \mathfrak{a}$ and assume that for each $%
X=X_{i}$, 
\begin{equation}
(m-1)\left( \dim \Phi -\dim \Psi \right) -1\geq m\left( \dim \Phi
_{X}-\min_{\sigma \in W}\dim (\Phi _{X}\cap \sigma (\Psi ))\right)
\label{Wr1}
\end{equation}%
for all $\mathbb{R}$-closed, root subsystems $\Psi \subseteq \Phi $ of
co-rank $1$. Then $\mu _{X_{1}}\ast \cdot \cdot \cdot \ast \mu _{X_{m}}$ is
absolutely continuous.
\end{corollary}

\begin{remark}
We note that if (\ref{Wr1}) holds for some $m$, then it holds for all
greater integers. We also note that if the inequality holds for $X$, then it
also holds for any $Y\in \mathfrak{a}$ with $\Phi _{Y}\subseteq \Phi _{X}$
since $\dim \Phi _{Y}-\dim (\Phi _{Y}\cap \Psi )\leq \dim \Phi _{X}-\dim
(\Phi _{X}\cap \Psi )$.
\end{remark}

The second result is a geometric characterization of absolute continuity. It
will be used both to prove absolute continuity (in this section) and
singularity in section 4.

\begin{theorem}
\label{keyprop} \cite{GHSym} Suppose $Z_{1},...,Z_{m}\in \mathfrak{a}$. Then 
$\mu _{Z_{1}}\ast \cdot \cdot \cdot \ast \mu _{Z_{m}}$ is absolutely
continuous if and only if there is some $k_{1}=Id$, $k_{2},...,k_{m}\in K$
(equivalently, for a.e. $(k_{1},...,k_{m})\in K^{m}$) such that 
\begin{equation}
sp\{Ad(k_{j})\mathcal{N}_{Z_{j}}:j=1,...,m\}=\mathfrak{p.}  \label{char}
\end{equation}
\end{theorem}

The nature of the proof of the sufficiency part of the theorem will differ
depending on the type of the exceptional symmetric spaces. There are three
cases to consider; these comprise subsections 3.2, 3.3, 3.4.

\subsection{The symmetric spaces whose restricted root spaces all have
dimension one [Types $EI$, $EV$, $EVIII$, $FI$ and $G$]}

We will first prove the absolute continuity of $\mu _{X}^{L(X)}$ for the
exceptional symmetric spaces in which all the restricted root spaces have
dimension one. These are the symmetric spaces in which the absolute and
restricted root systems coincide, the symmetric spaces $EI$, $EV$, $EVIII$, $%
FI$, and $G$. Their (absolute and restricted) root systems are types $E_{6}$%
, $E_{7}$, $E_{8}$, $F_{4}$, and $G_{2},$ respectively.

\begin{proposition}
\label{multone}Let $G/K$ be any of the symmetric spaces $EI$, $EV$, $EVIII$, 
$FI$, and $G$. Then inequality (\ref{Wr1}) holds with $m=L_{X}$, as defined
in Theorem \ref{main}, for any non-zero $X\in \mathfrak{a}$ and all $\Psi $
of co-rank $1.$
\end{proposition}

\begin{proof}
We note, first, that (\ref{Wr1}) is equivalent to the statement 
\begin{equation*}
(m-1)dim\mathcal{N}_{\Psi }-1\geq m(dim\Phi _{X}\cap \mathcal{N}_{\Psi })
\end{equation*}
for all $\mathbb{R}$-closed, root subsystems $\Psi \subseteq \Phi $ of
co-rank $1.$

Assume $\Phi $ has rank $n$ and base $\{\alpha _{1},..,\alpha _{n}\}$. Up to
Weyl conjugacy, a base for $\Psi $ consists of all but one $\alpha _{j}$.
Hence if $\{\lambda _{1},...,\lambda _{n}\}$ are the fundamental dominant
weights, then there exists $j$ such that $\Phi \backslash \sigma (\Psi
)=\{\alpha $ $\in \Phi :(\alpha ,\lambda _{j})\neq 0\}$ for a suitable Weyl
conjugate. Replacing $\Phi _{X}$ by $\Phi _{\sigma (X)}$, there is no loss
in assuming $\Phi \backslash \Psi =\{\alpha $ $\in \Phi :(\alpha ,\lambda
_{j})\neq 0\}$, thus our task is to prove%
\begin{equation}
(m-1)\dim \mathcal{N}_{\Psi _{j}}-1\geq m(\dim \Phi _{X^{\prime }}\cap 
\mathcal{N}_{\Psi _{j}})\text{ }  \label{redn}
\end{equation}%
whenever $\mathcal{N}_{\Psi _{j}}=sp\{E_{\alpha }^{-}:(\alpha ,\lambda
_{j})\neq 0\}$ and $X^{\prime }$ is Weyl conjugate to $X$.

This follows from the calculations done in \cite{HJSY} and \cite{HS}, as we
will now explain. Put 
\begin{equation*}
X_{1}=\{\alpha \in \Phi ^{+}:(\alpha ,\lambda _{j})\neq 0\}
\end{equation*}%
and 
\begin{equation*}
B_{1}=B_{1}(\Phi _{X^{\prime }})=\{\alpha \in \Phi _{X^{\prime
}}^{+}:(\alpha ,\lambda _{j})\neq 0\}.
\end{equation*}%
Thus $\dim \mathcal{N}_{\Psi _{j}}=\left\vert X_{1}\right\vert $ and $\dim
\Phi _{X^{\prime }}\cap \mathcal{N}_{\Psi _{j}}=\left\vert B_{1}\right\vert $%
, so (\ref{redn}) is equivalent to proving 
\begin{equation}
(m-1)\left\vert X_{1}\right\vert -m\left\vert B_{1}\right\vert \geq 1
\label{redn2}
\end{equation}%
for all $j$ and all $X^{\prime }$ Weyl conjugate to $X$.

It is shown in \cite[Section 2.2]{HS} that $m=2$ satisfies (\ref{redn2}),
for any \thinspace $j,$ when the root system is type $G_{2}$ or $F_{4},$ and
similarly $m=3$ works when the root system is type $E_{6}$, $E_{7}$ or $%
E_{8} $.

More refined calculations were done in \cite{HJSY} to show that the
inequality (\ref{redn2}) actually holds with $m=2,$ except when $\Phi _{X}$
is of type $D_{5}$ or $A_{5}$ in the root system $E_{6}$, type $E_{6}$ or $%
D_{6}$ in the root system $E_{7}$, or type $E_{7}$ in the root system $E_{8}$%
.
\end{proof}

To summarize:

\begin{corollary}
If $G/K$ is one of the symmetric spaces $EI$, $EV$, $EVIII$, $FI$, and $G$,
then $\mu _{X}^{L_{X}}$ is absolutely continuous for any non-zero $X\in 
\mathfrak{a}$. Furthermore, any convolution of two non-trivial orbital
measures in the symmetric spaces of type $FI$ or $G,$ or convolution of
three orbital measures in types $EI$, $EV$ and $EVIII$ is absolutely
continuous.
\end{corollary}

\subsection{An Embedding argument [Types $EII$, $EVI$ and $EIX$]}

Motivated by ideas in \cite{GSColloq}, it was shown in \cite{GHSym} that if $%
G_{1}/K_{1}$ suitably embeds into $G_{2}/K_{2}$, then absolute continuity
information for $G_{2}/K_{2}$ can be obtained from that for $G_{1}/K_{1}$.
In this section we will show that under mild technical conditions it is
possible to turn a suitable embedding of $(\mathfrak{g}_{1},\mathfrak{k}%
_{1}) $ into $(\mathfrak{g}_{2},\mathfrak{k}_{2})$ into an embedding of the
Lie groups' universal covers. This will be useful for us because (as we will
see) absolute continuity questions for orbital measures are unchanged if $G$
is replaced by its universal cover.

We will then apply this embedding argument to deduce the desired absolute
continuity properties for the symmetric space of Cartan types $EII$, $EVI$
and $EIX$ from those of the symmetric space $FI$ where we have already seen
in the previous subsection that $L_{X}=2$ holds.

We first recall what `suitably embeds' means at the group level and the
embedding theorem that was proven in \cite{GHSym}.

\begin{definition}
\label{groupembedding}Let $G_{1}/K_{1}$ and $G_{2}/K_{2}$ be two symmetric
spaces. We say that $G_{1}/K_{1}$ is \textbf{embedded} into $G_{2}/K_{2}$ if
there is a mapping $\mathcal{I}:G_{1}\rightarrow G_{2}$ satisfying the
following properties.

\begin{enumerate}
\item $\mathcal{I}$ is a group isomorphism into $G_{2}$.

\item $\mathcal{I}$ restricted to $A_{1}$ is a topological group isomorphism
onto $A_{2}$.

\item $\mathcal{I}$ maps $K_{1}$ into $K_{2}$.
\end{enumerate}
\end{definition}

\begin{theorem}
\label{embedding2}\cite{GHSym} Suppose $G_{1}/K_{1}$ is embedded into $%
G_{2}/K_{2}$ with the mapping $\mathcal{I}:G_{1}\rightarrow G_{2}$. Let $%
X_{1},....,X_{m}\in \mathfrak{a}_{1}$. If $\mu _{X_{1}}\ast \cdot \cdot
\cdot \ast \mu _{X_{m}}$ is absolutely continuous on $\mathfrak{p}_{1}$ and $%
\,\exp (Z_{j})=\mathcal{I}(\exp X_{j})$, then $\mu _{Z_{1}}\ast \cdot \cdot
\cdot \ast \mu _{Z_{m}}$ is absolutely continuous on $\mathfrak{p}_{2}$.
\end{theorem}

We now turn to embeddings of Lie algebras.

\begin{definition}
Let $\mathfrak{g}_{1}=\mathfrak{k}_{1}\oplus \mathfrak{p}_{1}$, $\mathfrak{g}%
_{2}=\mathfrak{k}_{2}\oplus \mathfrak{p}_{2}$ be non-compact, semi-simple
Lie algebras with the given Cartan decompositions and let $\mathfrak{a}%
_{i}\subseteq \mathfrak{p}_{i}$ be maximal abelian subalgebras. We say $(%
\mathfrak{g}_{1},\mathfrak{k}_{1})$ is \textbf{embedded} into $(\mathfrak{g}%
_{2},\mathfrak{k}_{2})$ if there is a Lie algebra homomorphism mapping $%
\iota :\mathfrak{g}_{1}\rightarrow \mathfrak{g}_{2}$ such that

\begin{enumerate}
\item $\iota$ is injective.

\item $\iota$ is an isomorphism from $\mathfrak{a}_1$ to $\mathfrak{a}_2$.

\item $\iota$ maps $\mathfrak{k}_1$ into $\mathfrak{k}_2$
\end{enumerate}
\end{definition}

It is well known that each Lie algebra is associated to exactly one simply
connected Lie group. We will now show that the above embedding descends to
the simply connected Lie groups associated to the Lie algebras. First, we
observe that the simply connected Lie group gives a symmetric space.

\begin{proposition}
Let $\mathfrak{g}$ be a non-compact, semi-simple Lie algebra with Cartan
involution $\theta $. Let $G$ be its associated simply connected, connected
Lie group and assume that $G$ has finite center $Z(G)$. Then there is an
involution $\Theta $ on $G$ with $d\Theta =\theta $ and if we let $K$ be the
Lie subgroup associated to $\mathfrak{k,}$ then $K$ is compact, connected,
contains $Z(G)$ and is the set of fixed points of $\Theta $.
\end{proposition}

\begin{proof}
This follows from ~\cite[VI Thm. 1.1]{He}. In particular, compactness of $K$
follows from the fact that $G$ has finite center.
\end{proof}

\begin{theorem}
\label{embedding1}Let $(\mathfrak{g}_{1},\mathfrak{k}_{1})$ be embedded into 
$(\mathfrak{g}_{2},\mathfrak{k}_{2})$ by $\iota :\mathfrak{g}_{1}\rightarrow 
\mathfrak{g}_{2}$. Suppose, in addition, that the simply connected Lie group 
$G_{1}$ associated to $\mathfrak{g}_{1}$ has finite center. For each Lie
group $G_{2},$ with Lie algebra $\mathfrak{g}_{2}$, there is a discrete
normal subgroup $N$ such that $G_{1}/N$ embeds into $G_{2}$.
\end{theorem}

\begin{proof}
Let $\theta _{i}$ be the Cartan involutions on $\mathfrak{g}_{i}$. Since $%
G_{1}$ is simply connected, $\theta _{1}$ descends to a unique involution $%
\Theta _{1}$ of $G_{1}$. Let $K_{1}$ be the Lie subgroup associated to $%
\mathfrak{k}_{1},$ the set of fixed points of $\Theta _{1}$ in $G_{1}$. Set $%
A_{1}=exp(\mathfrak{a}_{1})$.

Let $\mathcal{I}$ be the map from $G_{1}$ to $G_{2}$ such that $d\mathcal{I}%
_{e}=\iota $ (which exists because $G_{1}$ is simply connected). As $\iota $
is injective, $\mathcal{I}$ is locally injective, so its kernel is discrete
and therefore is contained in $Z(G_{1})\subseteq K$. As $Z(G_{1})$ is fixed
by $\Theta _{1},$ the map $\Theta _{1}$ induces an involution on any
quotient $G_{1}/N$ where $N\subseteq Z(G_{1})$.

Taking $N=\mathrm{ker}(\mathcal{I})$, gives the Cartan involution $\tilde{%
\Theta}_{1}$ on $G_{1}/N$ and the induced map $\tilde{\mathcal{I}}%
:G_{1}/N\rightarrow G_{2}$. Now $G_{1}/N$ still has Lie algebra $\mathfrak{g}%
_{1}$ and the image of $K_{1}$ under the quotient map $\pi $ has Lie algebra 
$\mathfrak{k}_{1}$. Also, note that $\pi (K_{1})\cong K_{1}/N$ and $\pi
(A_{1})=exp(\mathfrak{a}_{1})\cong A_{1}$ in $G_{1}/N$ since $K_{1}\cap N=N$
and $A_{1}\cap N=\{e\}$.

We claim $\tilde{\mathcal{I}}$ embeds $G_{1}/N$ into $G_{2}$ in the sense of
Definition \ref{groupembedding}.

(1) $\tilde{\mathcal{I}}$ is a group homomorphism into $G_{2}$ and clearly
it is injective.

(2) Since $K_{i}$ is generated by $exp(\mathfrak{k}_{i})$ and $\iota $ maps $%
\mathfrak{k}_{1}$ into $\mathfrak{k}_{2}$, $\mathcal{I}$ maps $K_{1}$ into $%
K_{2}$ and $\tilde{\mathcal{I}}$ maps $\pi (K_{1})$ into $K_{2}$.

(3) Similarly, $\mathcal{I}$ maps $\pi (A_{1})\cong A_{1}$ into $A_{2}$. It
is surjective since for any $exp(Y)\in A_{2}$, $\mathcal{I}(exp(\iota
^{-1}(Y)))=exp(Y)$ and hence is a bijection from $\pi (A_{1})\cong A_{1}$ to 
$A_{2}$. Because $\varphi :\mathfrak{p}_{i}\times K_{i}\rightarrow G_{i}$
given by $\varphi (X,k)=(\exp X)k$ is a diffeomorphism (\cite[VI 1.1]{He}),
it is an open map and thus $A_{i}=\varphi (\mathfrak{a}_{i},{e})$ is closed
being the image of a closed set. As $\pi $ is an open map, $\pi (A_{1})$ is
also closed.

Thus $\pi (A_{1})$ and $A_{2}$ are Lie groups with the subspace topology and 
$\tilde{\mathcal{I}}$ is a bijective Lie group homomorphism between the two.
Hence $\mathcal{I}:\pi (A_{1})\rightarrow A_{2}$ is a bijective local
diffeomorphism, which means it is a global diffeomorphism, and so a
topological group isomorphism.
\end{proof}

To apply this result, we will make use of the Freudenthal magic square
construction for models of the exceptional Lie algebras/groups. The
construction we will use is due to Vinberg \cite{OV}. We give a brief
overview here, as explained in Barton and Sudbery \cite{BS}.

The method we use is called the Cayley-Dickson construction. We begin with $F%
\mathbb{=R},\mathbb{C}$ or the quaternions $\mathbb{H}$, and define a
multiplication on $F\times F$ by $(x_{1},y_{1})(x_{2},y_{2})=(x_{1}x_{2}+%
\varepsilon \overline{y_{2}}y_{1},y_{2}x_{1}+y_{1}\overline{x_{2}})$, with
conjugation given by $\overline{(x,y)}=(\overline{x},-y)$ and quadratic form 
$N(x,y)=(x,y)\overline{(x,y)}$. Here $\varepsilon =\pm 1$, with $\varepsilon
=-1$ giving a division algebra and $\varepsilon =1$ giving a split algebra.

We will write $F\oplus lF$ to denote the composition algebra that arises
from the Cayley-Dickson construction starting with $F$ and the new imaginary
unit $l$ where $l^{2}=\varepsilon $. Of course, when $\varepsilon =-1$ and $%
F=\mathbb{R}$, we have $l=i$ and $F\oplus lF=\mathbb{C}$. We can similarly
obtain $\mathbb{H}$ and the octonions, $\mathbb{O}$, in this fashion.

Given unital composition algebras $C_{j}$ over $\mathbb{R}$, we let $%
A_{3}(C_{1}\otimes C_{2})$ denote the $3\times 3$ matrices with entries in $%
C_{1}\otimes C_{2}$ that are skew-Hermitian and let $A_{3}^{\prime
}(C_{1}\otimes C_{2})$ denote the subset of trace zero matrices. We denote
by $Der(C_{j})$ the derivations on $C_{j}$.

Put%
\begin{equation}
V_{3}(C_{1},C_{2})=A_{3}^{\prime }(C_{1}\otimes C_{2})\oplus
Der(C_{1})\oplus Der(C_{2}).  \label{Vinberg}
\end{equation}

With the definitions given below for the Lie bracket, $V_{3}(C_{1},C_{2})$
is a Lie algebra.

(1) $Der(C_{1})\oplus Der(C_{2})$ is a Lie subalgebra with $%
[Der(C_{1}),Der(C_{2})]=0$ and the Lie bracket within each subalgebra is the
standard Lie bracket, $[X,Y]=XY-YX$.

(2) For $D\in Der(C_{1})\oplus Der(C_{2})$ and $A\in A_{3}^{\prime
}(C_{1}\otimes C_{2})$ we put $[D,A]=D(A),$ where $D(A)$ is defined to be
the matrix whose entries are obtained by applying $D$ to the corresponding
entry in $A$.

(3) For $A,B\in A_{3}^{\prime }(C_{1}\otimes C_{2})$ we have 
\begin{equation*}
\lbrack A,B]=AB-BA-\frac{tr(AB-BA)}{3}I+\frac{1}{3}%
\sum_{i,j}D_{A_{ij},B_{ij}}
\end{equation*}%
where 
\begin{equation*}
D_{u_{1}\otimes v_{1},u_{2}\otimes v_{2}}=\langle v_{1},v_{2}\rangle
D_{u_{1},u_{2}}+\langle u_{1},u_{2}\rangle D_{v_{1},v_{2}}
\end{equation*}%
and 
\begin{equation*}
D_{x,y}=[L_{x},L_{y}]+[L_{x},R_{y}]+[R_{x},R_{y}].
\end{equation*}%
Here $L$ and $R$ denote left and right multiplication and the brackets are
the standard brackets for operators on a vector space.

This Lie algebra is relevant for us because with appropriate choices of $%
C_{1}$ and $C_{2}$ it gives rise to various exceptional symmetric spaces. In
particular, if $C_{1}=\mathbb{H\oplus }$ $l\mathbb{H}$ is the split
octonions and $C_{2}$ is a division algebra, then we obtain the exceptional
symmetric spaces as indicated below:

\medskip

\begin{center}
\begin{tabular}{|p{2.5cm}|c|c|c|c|}
\hline
Choice for $C_2$ & $\mathbb{R}$ & $\mathbb{C}$ & $\mathbb{H}$ & $\mathbb{O}$
\\ \hline
Cartan type of $V_3(C_1,C_2)$ & $FI$ & $EII$ & $EVI$ & $EIX$ \\ \hline
\end{tabular}
\end{center}

\medskip

We will need to understand the Cartan decomposition of these algebras.
Towards this end, consider a composition algebra $C_{1}=F\oplus lF$. By $%
F^{\prime }$ we mean the subspace of $F$ orthogonal to $\mathbb{R}$.

For each linear, symmetric (and traceless if $F=\mathbb{H}$) map $%
S:F^{\prime }\rightarrow F^{\prime },$ define $G_{S}$ by 
\begin{equation*}
G_{S}(a)=l(Sa)\text{ for }a\in F^{\prime }\text{ and }G_{S}(l)=0,
\end{equation*}
and extend $G_{S}$ to $C_{1}$ by linearity and the derivation property. It
can be checked that $G_{S}$ is a derivation on $C_{1}$ and that if $[\cdot
,\cdot ]$ denotes the usual bracket on $Der(C_{1})$, then $%
[G_{S},G_{T}]=l^{2}\overline{[S,T]},$ where we denote by $\overline{D}$ the
unique extension of $D$ to a derivation on $C_{1}$ that sends $l$ to $0$.

The derivations $\,E_{a}$ and $F_{a}$ mentioned in the next theorem are
defined in \cite{BS} for each $a\in F^{\prime }$, but as the specific
details of these will not be important for us we have not given the
definitions here.

\begin{theorem}
\cite[Section 7]{BS} Let $C_{1}=F\oplus lF$ be a split algebra and $C_{2}$
be a division algebra. Let $V_{3}(C_{1},C_{2})$ denote the Lie algebra as in
(\ref{Vinberg}). The Cartan decomposition $V_{3}(C_{1},C_{2})=\mathfrak{k}%
\oplus \mathfrak{p}$ is given by 
\begin{equation*}
\mathfrak{k}=A_{3}^{\prime }(F\otimes C_{2})\oplus Der_{0}(C_{1})\oplus
Der(C_{2})
\end{equation*}%
and 
\begin{equation*}
\mathfrak{p}=A_{3}^{\prime }(lF\otimes C_{2})\oplus Der_{1}(C_{1}),
\end{equation*}%
where%
\begin{equation*}
Der_{0}(C_{1})=\text{sp}\{\overline{D},E_{a}:D\in Der(F),a\in F^{\prime }\}%
\text{ }
\end{equation*}%
and%
\begin{equation*}
Der_{1}(C_{1})=\text{sp}\{F_{a},G_{S}:a\in F^{\prime },S\text{ linear,
symmetric}\}.
\end{equation*}
\end{theorem}

Using this decomposition, we can find a maximal abelian subalgebra $%
\mathfrak{a}$ in $\mathfrak{p}$. This is the content of the next proposition.

\begin{proposition}
\label{a} Let $V_{3}(C_{1},C_{2})=\mathfrak{k}\oplus \mathfrak{p}$ be a Lie
algebra obtained from the Vinberg construction, with Cartan decomposition as
given above and $C_{1}=\mathbb{H}\oplus l\mathbb{H}$. Define the following
elements in $A_{3}^{\prime }(l\mathbb{H}\otimes C_{2})$: 
\begin{equation*}
H_{1}=%
\begin{bmatrix}
l\otimes 1 & 0 & 0 \\ 
0 & -l\otimes 1 & 0 \\ 
0 & 0 & 0%
\end{bmatrix}%
,\text{ }H_{2}=%
\begin{bmatrix}
0 & 0 & 0 \\ 
0 & l\otimes 1 & 0 \\ 
0 & 0 & -l\otimes 1%
\end{bmatrix}%
.
\end{equation*}%
A maximal abelian subalgebra is spanned by $H_{1}$, $H_{2}$ and $G_{S_{1}}$, 
$G_{S_{2}},$ where $S_{1},S_{2}$ span the diagonal traceless operators $%
\mathbb{H}^{\prime }\rightarrow \mathbb{H}^{\prime }$.
\end{proposition}

\begin{proof}
First, observe that $D_{1,1}=0$ since $L_{1}$, $R_{1}$ are both commutative
and $D_{l,l}x=[L_{l},R_{l}]x=l(xl)-(lx)l=0$. It is clear that $%
H_{1}H_{2}-H_{2}H_{1}=0$ so we have that $[H_{1},H_{2}]=0$.

As $G_{S_{1}}l=0$ we have 
\begin{equation*}
\lbrack G_{S_{1}},H_{j}]=G_{S_{1}}(H_{j})=0\text{ for }j=1,2.
\end{equation*}
Also, $[G_{S_{1}},G_{S_{2}}]=l^{2}\overline{[S_{1},S_{2}]}=0$ since both $%
S_{1}$, $S_{2}$ were chosen to be diagonal. Thus the given vectors span an
abelian subalgebra.

Maximality follows because the maximal abelian subalgebra has dimension
equal to the rank of its associated symmetric space.
\end{proof}

Having explained the Vinberg construction, we will now show there are
embeddings of the Lie algebras and then appeal to Theorem \ref{embedding1}
to obtain suitable embeddings at the group level. This could be done in more
generality, but it is the case $F=\mathbb{H}$ that is of interest to us.

\begin{proposition}
Let $C_{0}=\mathbb{H}\oplus l\mathbb{H}$ be the split octonions and let $%
C_{1}$ denote one of $\mathbb{C}$, $\mathbb{H}$ or $\mathbb{O}$. Then $%
V_{3}(C_{0},\mathbb{R})$ embeds into $V_{3}(C_{0},C_{1}).$
\end{proposition}

\begin{proof}
Let $\mathfrak{t}_{0}$ $\oplus $ $\mathfrak{p}_{0}$ and $\mathfrak{t}_{1}$ $%
\oplus $ $\mathfrak{p}_{1}$ denote the Cartan decompositions of $V_{3}(C_{0},%
\mathbb{R})$ and $V_{3}(C_{0},C_{1})$ respectively, and let $\mathfrak{a}%
_{0} $ and $\mathfrak{a}_{1}$ be associated the maximal abelian subalgebras
described in Prop. \ref{a}. We note that $V_{3}(C_{0},\mathbb{R}%
)=A_{3}^{\prime }(C_{0}\otimes \mathbb{R})\oplus Der(C_{0})$ as $\mathbb{R}$
has no non-trivial derivations.

Define the map $\iota :V_{3}(C_{0},\mathbb{R})\rightarrow V_{3}(C_{0},C_{1})$
by 
\begin{equation*}
\iota (A+D)=A+D.
\end{equation*}
Here we have the natural inclusion $A_{3}^{\prime }(C_{0}\otimes \mathbb{R}%
)\rightarrow A_{3}^{\prime }(C_{0}\otimes C_{1})$, which is just an
extension of scalars, and the identity map on $Der(C_{0})$.

One can easily check that $\iota $ is a Lie algebra homomorphism.
Furthermore, $\iota $ is clearly an injective map into $V_{3}(C_{0},C_{1})$.
Since $\iota (H_{i})=H_{i}$ and $\iota (G_{S_{i}})=G_{S_{i}}$ (with the
natural identifications), it is clear from Prop. \ref{a} that $\iota :%
\mathfrak{a}_{0}\rightarrow \mathfrak{a}_{1}$ bijectively. Finally, $\iota $
maps $A_{3}^{\prime }(\mathbb{H}\otimes \mathbb{R})$ into $A_{3}^{\prime }(%
\mathbb{H}\otimes C_{1})$, $Der_{0}(C_{0})$ to itself, hence it maps $%
\mathfrak{k}_{0}$ into $\mathfrak{k}_{1}$.
\end{proof}

\begin{corollary}
We have the embeddings of the symmetric space of Cartan type $FI$ into $EII$%
, $EVI$ and $EIX.$
\end{corollary}

\begin{proof}
Let $G$ be the non-compact simple Lie group with Lie algebra $V_{3}(C_{0},%
\mathbb{R})$ and let $\widetilde{G}$ denote its simply connected universal
cover. The centre of $\widetilde{G}$ is the fundamental group of $G$, namely 
$\mathbb{Z}_{2}$ \cite{Car}, thus Theorem \ref{embedding1} lifts the
embeddings of the Vinberg construction demonstrated in the previous
proposition to the group level.
\end{proof}

\begin{corollary}
The convolution of any two orbital measures in $EII$, $EVI$ or $EIX$ is
absolutely continuous.
\end{corollary}

\begin{proof}
We first remark that the problem of the absolute continuity of the product
of orbital measures is unaffected by the choice of $G$ with a given Lie
algebra $\mathfrak{g}$. This is because any such group $G$ is isomorphic to $%
\widetilde{G}/Z$ where $\widetilde{G}$ is the common universal cover and $Z$
is a subset of the centre of $\widetilde{G}.$ Further, $K\simeq \widetilde{K}%
/Z.$ As the adjoint action is trivial on the centre, the adjoint action on $%
K $ and $\widetilde{K}$ coincide.

As shown in Prop. \ref{multone}, the convolution of any two orbital measures
in $FI$ is absolute continuous. Thus the corollary holds because of the
embedding Theorem \ref{embedding1}.
\end{proof}

\subsection{The remaining exceptional symmetric spaces [Types $EIII$, $EIV$, 
$EVII$, $FII$]}

For the remaining exceptional symmetric spaces we use ad-hoc arguments to
prove the required absolute continuity results.

Type $FII$: Here the restricted root space is type $BC_{1}$. Any non-zero $%
X\in \mathfrak{a}$ is regular, meaning $\Phi _{X}$ is empty, and it is known
that the convolution of any two regular elements is absolutely continuous 
\cite{GS2003}.

Type $EIV$: The restricted root space is type $A_{2}$. All non-zero $X$ are
regular, except $X$ of type $A_{1}$, and for the regular $X$, $\mu _{X}^{2}$
is absolutely continuous. As the rank of $EIV$ is $2,$ any convolution
product of three non-trivial orbital measures is absolutely continuous
according to \cite{GSJFA}.

Type $EIII$: The restricted root space is type $BC_{2},$ with the roots $%
\alpha _{1}\pm \alpha _{2}$ having multiplicity $6$, the roots $\alpha _{j}$
having multiplicity $8$ and the roots $2\alpha _{j}$ having multiplicity $1$%
. The co-rank one (equivalently, rank one in this case) closed root
subsystems are either $\{\alpha _{1}+\alpha _{2}\}$, $\{\alpha _{1}-\alpha
_{2}\}$ or $\{\alpha _{j},2\alpha _{j}\}$. With this information it is
simple to check that (\ref{Wr1}) is satisfied with $m=2$ for all $\Phi _{X},$
$X\neq 0$. Consequently any product of two $K$-invariant, continuous orbital
measures is absolutely continuous.

Type $EVII$: In this case, the restricted root system is type $C_{3},$ with
the short roots having multiplicity $8$, the long roots multiplicity $1$,
and $\dim \Phi =51$. The co-rank one, closed root subsystems $\Psi $ are
type $C_{2}$, $A_{2}$ or $A_{1}\times C_{1}$.

If $X$ is of type $C_{2}$, then $\dim \Phi _{X}=18$. The intersection of any
such $\Phi _{X}$ with any $\Psi $ of type $A_{2}$ must contain a short root
and hence the intersection has dimension at least $8$. It easily follows
from (\ref{Wr1}) that $\mu _{X}^{3}$ is absolutely continuous. If $X$ is of
type $A_{1}\times C_{1}$, then $\dim \Phi _{X}=9$ and it is even easier to
see that in this case $\mu _{X}^{2}$ is absolutely continuous. If $X$ is
type $A_{1}$ or $C_{1}$, then $\Phi _{X}$ is contained in the set of
annihilating roots of an $X^{\prime }$ of type $A_{1}\times C_{1}$ and again
it follows that $\mu _{X}^{2}$ is absolutely continuous.

Finally, if $X$ is type $A_{2}$, then $\dim \Phi _{X}=24$ and (\ref{Wr1}) is
clearly satisfied with $m=3$. It follows that the convolution of any $3$
non-trivial orbital measures on $\mathfrak{p}$ is absolutely continuous.
Thus $L(G)=3$ suffices.

However, to complete the proof of the absolute continuity claims of Theorem %
\ref{main} we must show that if $X$ is type $A_{2}$, then $\mu _{X}^{2}$ is
absolutely continuous. For this we will use the geometric criteria for
absolute continuity given in Theorem \ref{keyprop}.

The absolute root system of the symmetric space $EVII$ is type $E_{7}$ and
we can take as a base the roots $\frac{1}{2}(e_{1}+e_{8}-(e_{2}-\cdot \cdot
\cdot -e_{7})$, $e_{1}+e_{2}$, $e_{1}-e_{2}$, $e_{3}-e_{2}$, $%
...,e_{6}-e_{5} $. One can see from the Satake diagram (c.f., \cite[p. 534]%
{He}) that $e_{j}|_{\mathfrak{a}}=0$ for $j=1,2,3,4$ and that the roots $%
e_{5}\pm e_{6}$ and $e_{7}-e_{8}$ are distinct non-zero elements of $%
\mathfrak{a}^{\ast }$. We have $\theta (e_{j})=e_{j}$ for $j=1,2,3,4,$%
\begin{equation*}
\theta (e_{5}\pm e_{6})=-(e_{5}\pm e_{6}),\text{ }\theta
(e_{7}-e_{8})=-(e_{7}-e_{8}),
\end{equation*}%
\begin{equation*}
\theta (e_{j}\pm e_{i})=-e_{j}\pm e_{i}\text{ for }j=5,6\text{, }i=1,2,3,4
\end{equation*}%
and%
\begin{equation*}
\theta \left( \frac{1}{2}(e_{8}-e_{7}+e_{6}\pm
e_{5}+\sum_{j=1}^{4}s_{j}e_{j})\right) =-\frac{1}{2}\left(
e_{8}-e_{7}+e_{6}\pm e_{5}-\sum_{j=1}^{4}s_{j}e_{j}\right)
\end{equation*}%
(here $s_{j}=\pm 1$ and $\prod s_{j}=-1$).

Up to Weyl conjugacy, the set of annihilating roots of an element $X$ of
type $A_{2}$ consists of the restricted roots $\{\frac{1}{2}%
(e_{8}-e_{7}+e_{6}\pm e_{5}),$ $e_{5}\}$, each of which has multiplicity $8$%
. The non-annihilating roots are the restricted roots $\{\frac{1}{2}%
(e_{8}-e_{7}-e_{6}\pm e_{5}),$ $e_{6},$ $e_{5}\pm e_{6},$ $e_{8}-e_{7}\},$
with the latter three having multiplicity one. The restricted root spaces of
the non-annihilating roots of multiplicity 8 are%
\begin{eqnarray*}
\mathfrak{g}_{\frac{1}{2}(e_{8}-e_{7}-e_{6}\pm e_{5})} &=&sp\left\{ E_{\beta
}^{-}:\beta =\frac{1}{2}(e_{8}-e_{7}-e_{6}\pm
e_{5}+\sum_{j=1}^{4}s_{j}e_{j})\right\} \text{ and} \\
\mathfrak{g}_{e_{6}} &=&sp\{E_{e_{6}\pm e_{j}}^{-}:j=1,...,4\}.
\end{eqnarray*}

Let 
\begin{eqnarray*}
F_{1} &=&E_{e_{5}+e_{6}}+E_{-(e_{5}+e_{6})}=E_{e_{5}+e_{6}}^{+}, \\
F_{2} &=&E_{e_{5}-e_{6}}+E_{e_{6}-e_{5}}=E_{e_{5}-e_{6}}^{+}\text{,} \\
F_{3} &=&E_{e_{8}-e_{7}}+E_{e_{7}-e_{8}}=E_{e_{8}-e_{7}}^{+}.
\end{eqnarray*}%
Routine calculations show that if $F=c_{1}F_{1}+c_{2}F_{2}+c_{3}F_{3}\in K$,
then 
\begin{equation*}
ad(F)E_{e_{6}\pm e_{j}}^{-}=(c_{2}-c_{1})E_{e_{5}\pm e_{j}}^{-}\text{ for }%
j=1,2,3,4.
\end{equation*}%
If 
\begin{equation*}
\beta ^{+}=\frac{1}{2}(e_{8}-e_{7}+e_{6}-e_{5}+\sum_{j=1}^{4}s_{j}e_{j}),%
\text{ }\beta ^{-}=\frac{1}{2}(e_{8}-e_{7}-e_{6}+e_{5}+%
\sum_{j=1}^{4}s_{j}e_{j}),
\end{equation*}%
and 
\begin{equation*}
\gamma ^{+}=\frac{1}{2}(e_{8}-e_{7}+e_{6}+e_{5}+\sum_{j=1}^{4}s_{j}e_{j}),%
\text{ }\gamma ^{-}=\frac{1}{2}(e_{8}-e_{7}-e_{6}-e_{5}+%
\sum_{j=1}^{4}s_{j}e_{j}),
\end{equation*}%
then $\beta ^{-}|_{\mathfrak{a}},\gamma ^{-}|_{\mathfrak{a}}$ are
non-annihilating roots and%
\begin{equation*}
ad(F)E_{\beta ^{-}}^{-}=(c_{2}-c_{3})E_{\beta ^{+}}^{-}\text{ and }%
ad(F)E_{\gamma ^{-}}^{-}=(c_{1}-c_{3})E_{\gamma ^{+}}^{-}\text{ .}
\end{equation*}%
Furthermore, if $c_{j}\neq 0$, then 
\begin{equation*}
sp(ad(F)\{E_{e_{5}\pm e_{6}}^{-},E_{e_{8}-e_{7}}^{-}\})=\mathfrak{a.}
\end{equation*}%
Thus, provided we choose the $c_{j}$ non-zero and distinct, then 
\begin{equation*}
\mathcal{N}_{X}+ad(F)\mathcal{N}_{X}=\mathfrak{p}\text{.}
\end{equation*}

Now consider $f_{t}=\exp tF\in K$ for small $t>0$. For any $Y$, $%
Ad(f_{t})(Y)=Y+t\cdot ad(F)(Y)+S_{t}(Y)$ where $\left\Vert S_{t}\right\Vert
\leq Ct^{2}$ for a constant $C$ independent of $Y$ and $t$. Hence, for all $%
t>0,$%
\begin{equation*}
sp\{\mathcal{N}_{X},Ad(f_{t})\mathcal{N}_{X}\}=sp\left\{ \mathcal{N}%
_{X},(ad(F)+\frac{1}{t}S_{t})(\mathcal{N}_{X})\right\} \subseteq \mathfrak{p}
\end{equation*}%
where $\left\Vert \frac{1}{t}S_{t}\right\Vert \leq Ct$. But since $sp\{%
\mathcal{N}_{X},ad(F)\mathcal{N}_{X}\}=\mathfrak{p}$, a linear algebra
argument implies that the same is true for $sp\{\mathcal{N}_{X},Ad(f_{t})%
\mathcal{N}_{X}\}$ for small enough $t>0$. Thus Theorem \ref{keyprop}
implies $\mu _{X}^{2}$ is absolutely continuous when $X$ is type $A_{2}$.

That completes the absolute continuity part of Theorem \ref{main}.

\section{Singularity}

We now turn to proving that $\mu _{X}^{2}$ is singular to Lebesgue measure
in the cases where we have claimed $L_{X}=3$. (Of course, if $L_{X}=2$, then
that is sharp as all orbits have measure zero and hence all orbital measures
are singular.)

When $X$ is type $D_{5}$ in $EI$, type $E_{6}$ or $D_{6}$ in $EV$, or type $%
E_{7}$ in $EVIII$, this is very easy: Simply note that for all of these
elements $X$, $\dim \mathcal{N}_{X}<\frac{1}{2}\dim \mathfrak{p}$ and hence
inequality (\ref{char}) of Theorem \ref{keyprop}, with $m=2,$ must
necessarily fail.

For the remaining cases, we will make use of the following condition that
guarantees singularity. It was originally established in \cite{Wr} for
compact Lie algebras.

\begin{proposition}
Let $X_{1},...,X_{m}\in \mathfrak{a}$ and assume that for $j=1,...,m$, $%
\sigma _{j}(\mathcal{N}_{X_{j}})\cap \mathcal{N}_{\Psi }$, are pairwise
disjoint for some closed, co-rank one root subsystem $\Psi $ of $\Phi $ and $%
\sigma _{j}\in W$. Then $\mu _{X_{1}}\ast \cdot \cdot \cdot \ast \mu
_{X_{m}} $ is singular to Lebesgue measure.
\end{proposition}

\begin{proof}
The proof in the symmetric space case is very similar to that given in \cite%
{Wr}. Given $\Psi ,$ choose $Z\in \mathfrak{a}$ satisfying $\Phi _{Z}=\Psi $%
. For each $j$, let $K_{j}=\{k\in K:Ad(k)X_{j}=X_{j}\}$ and let $\mathfrak{p}%
_{j}=\{Y\in \mathfrak{p:}[Y,X_{j}]=0\}$. Define $M_{j}$ to be the $K_{j}$%
-orbit of $Z$. Then $M_{j}$ is a submanifold of $O_{Z}$ and is contained in $%
\mathfrak{p}_{j}$. The tangent space to $M_{j}$ at $Z$ is 
\begin{equation*}
T_{Z}(M_{j})=\{Y\in \mathfrak{p:}[Y,Z]\neq 0\}\cap \mathfrak{p}_{j}=\mathcal{%
N}_{Z}\cap \mathfrak{p}_{j}\text{.}
\end{equation*}%
Note that $(T_{Z}(M_{j}))^{\bot }$ in $T_{Z}(O_{Z})$ is equal to $\mathcal{N}%
_{Z}\cap \mathcal{N}_{X_{j}}$. By \cite[Lemma 2.5]{Wr}, $\bigcap%
\limits_{j=1}^{m}Ad(k_{j})M_{j}$ is non-empty for $(k_{1},...,k_{m})$ near
the identity and hence the same is true for the larger space $%
\bigcap\limits_{j=1}^{m}Ad(k_{j})\mathfrak{p}_{j}$. But then standard
Hilbert space arguments (c.f., \cite[Lemma 4]{GHSym}) imply that 
\begin{equation*}
sp\{Ad(k_{j})\mathcal{N}_{X_{j}}:j=1,...,m\}\neq \mathfrak{p}
\end{equation*}%
for all $(k_{1},...,k_{m})$ near the identity. By Theorem \ref{keyprop} this
implies $\mu _{X_{1}}\ast \cdot \cdot \cdot \ast \mu _{X_{m}}$ is singular.
\end{proof}

We will use this criteria for the remaining cases.

\begin{enumerate}
\item $X$ type $A_{1}$ in $EIV$ where $\Phi $ is type $A_{2}$: By passing to
a Weyl conjugate, we can assume $\Phi _{X}=\{e_{1}-e_{2}\}$ and $\Phi
_{\sigma (X)}=\{e_{2}-e_{3}\}$. Put $\Psi =\{e_{1}-e_{3}\}$. Then $\mathcal{N%
}_{X}\cap \mathcal{N}_{\Psi }=\{e_{1}-e_{3}\},$ while $\mathcal{N}_{\sigma
(X)}\cap \mathcal{N}_{\Psi }=\{e_{1}-e_{2}\}$.

\item $X$ type $C_{2}$ in $EVII$ where $\Phi $ is type $C_{3}$: Here we can
assume $\Phi _{X}=\{e_{1}\pm e_{2},2e_{1},2e_{2}\}$, $\Phi _{\sigma
(X)}=\{e_{2}\pm e_{3},2e_{2},2e_{3}\},$ and take $\Psi =\{e_{1}\pm
e_{3},2e_{1},2e_{3}\}.$

\item $X$ type $A_{5}$ in $EI$ where $\Phi $ is type $E_{6}$: This was shown
to be singular using the criterion of Theorem \ref{keyprop} for the Lie
algebra setting in \cite[4.3]{HJSY}. The same argument applies here.
\end{enumerate}

This completes the singularity argument and hence the proof of Theorem \ref%
{main}.

\begin{landscape}
\centering
\section*{Basic facts about non-compact irreducible Riemannian symmetric spaces}

In the columm labelled $\mf{g}$, the number in the backet is the signature of the Killing form, dim $\mf{p} -$ dim $\mf{k}$.
\medskip

 \begin{tabular}{|c||c| c| >{\centering}m{2.25cm}| >{\centering}m{2.25cm}| >{\centering}m{1.75cm}| >{\centering}m{1.5cm}| c|} 
 \hline
  & $\mf{g}$ & $\mf{k}$ & Absolute Root System & Restricted Root System & Dimension of $G/K$ & Rank of $G/K$ & Multiplicities  \\ [0.5ex] 
 \hline\hline
 $EI$ & $\mf{e}_{6(6)}$ & $\mf{sp}(4)$ & $E_6$ &$E_6$ & 42 & 6 & All: 1  \\ [5pt]
 \hline
  $EII$ & $\mf{e}_{6(2)}$ & $\mf{su}(6)\oplus \mf{su}(2)$ & $E_6$ &$F_4$ & 40 & 4 & Not Needed \\ [5pt]
 \hline
  $EIII$ & $\mf{e}_{6(-14)}$ & $\mf{so}(10)\oplus \mf{so}(2)$ & $E_6$ &$BC_2$ & 32 & 2 & \makecell{$\varepsilon _i\pm\varepsilon _j$: 6 \\ $\varepsilon _i$: 8 \\$2\varepsilon _i$: 1} \\ [5pt]
 \hline
  $EIV$ & $\mf{e}_{6(-26)}$ & $\mf{f}(4)$ & $E_6$ &$A_2$ & 26 & 2 &  All: 8\\ [5pt]
 \hline
  $EV$ & $\mf{e}_{7(7)}$ & $\mf{su}(8)$ & $E_7$ &$E_7$ & 70 & 7 &  All: 1\\ [5pt]
 \hline
  $EVI$ & $\mf{e}_{7(-5)}$ & $\mf{so}(12)\oplus \mf{su}(2)$ & $E_7$ &$F_4$ & 64 & 4 & Not Needed \\ [5pt]
 \hline
  $EVII$ & $\mf{e}_{7(-25)}$ & $\mf{e}(6)\oplus \mf{so}(2)$ & $E_7$ &$C_3$ & 54 & 3 & \makecell{$\varepsilon _i \pm \varepsilon _j$: 8 \\ $2\varepsilon _i$: 1} \\ [5pt]
 \hline
  $EVIII$ & $\mf{e}_{8(8)}$ & $\mf{so}(16)$ & $E_8$ &$E_8$ & 128 & 8 & All: 1 \\ [5pt]
 \hline
  $EIX$ & $\mf{e}_{8(-24)}$ & $\mf{e}(7)\oplus \mf{su}(2)$ & $E_8$ &$F_4$ & 112 & 4 & Not Needed \\ [5pt]
 \hline
  $FI$ & $\mf{f}_{4(4)}$ & $\mf{sp}(3)\oplus \mf{su}(2)$ & $F_4$ &$F_4$ & 28 & 4 & All: 1 \\ [5pt]
 \hline
  $FII$ & $\mf{f}_{4(-20)}$ & $\mf{so}(9)$ & $F_4$ &$BC_1$ & 16 & 1 & \makecell{$\varepsilon _1$: 8\\ $2\varepsilon _1$: 7} \\ [5pt]
 \hline
 $G$ & $\mf{g}_{2(2)}$ & $\mf{su}(2)\oplus\mf{su}(2)$ & $G_2$ &$G_2$ & 8 & 2 & All: 1 \\ [5pt]
 \hline

\end{tabular}
\end{landscape}

\end{document}